\documentstyle[12pt]{article}

  \newcommand{\const}{\rm const}
  \newcommand{\Var}{\rm Var}
  \newcommand{\supp}{\rm supp}
  
  \newcommand{\vraisup}{\rm vraisup}
  \newcommand{\Cov}{\rm Cov}
  \newcommand{\Dom}{\rm Dom}
  
  \newcommand{\argmax}{\rm argmax}

\begin{document}

   \begin{center}

  \   {\bf   Covariation inequality  in  Grand Lebesgue Spaces. }\\

\vspace{4mm}

  \   {\bf Ostrovsky E., Sirota L.}\\

\vspace{4mm}

 ISRAEL,  Bar - Ilan University, department  of Mathematic and Statistics, 59200, \\

\vspace{4mm}

E - mails: eugostrovsky@list.ru, \hspace{5mm} sirota3@bezeqint.net \\

\vspace{5mm}

  {\bf Abstract} \\

\vspace{4mm}

 \end{center}

 \  We represent  in this preprint  the  exact estimate for covariation berween two random variables (r.v.),
which are measurable relative  the corresponding sigma - algebras through anyhow mixing coefficients. \par
 \  We associate a solution of this problem with fundamental function for correspondent rearrangement
invariant spaces. \par

\vspace{4mm}


\vspace{5mm}

\section{ \ Definitions.  Notations. Previous results.  Statement of problem.}

\vspace{4mm}

 \ Let  $ \  (\Omega, {\cal{B}}, {\bf P})  \  $ be probability space with correspondent expectation $ \ {\bf E}, \ $
variance $ \  \Var \ $ and covariation  $ \Cov: \ $

$$
\Cov(\xi, \eta) := {\bf E} \xi \eta - {\bf E \xi} \ {\bf E \eta}.
$$

\vspace{4mm}

 \  {\bf Denotation 1.1.} \par

{\it  \ We denote for arbitrary  sub - sigma algebra (field)
$ \   F \subset  {\cal B} \  $ and for arbitrary numerical valued random variable $ \ \xi \ $ the symbol}

$$
\xi \in \in F \eqno(1.0)
$$
{\it iff the r.v. $ \ \xi \ $ is measurable relative the sigma - field } $ \ F. \ $ \par

\vspace{4mm}

 \ Let $  \ F \  $ and $ \ G \ $ be two  sub - sigma algebras of source sigma field $ \ {\cal B}. \ $  We define as ordinary
the  so - called  uniform mixing coefficient, or equally Rosenblatt's coefficient $ \  \alpha(F,G) \  $ by the formula

$$
\alpha = \alpha(F,G)  \stackrel{def}{=} \sup_{ A \in F, \  B \in G }  |{\bf P}(A B) - {\bf P}(A) \ {\bf P}(B) |. \eqno(1.1)
$$

 \ The strong  mixing coefficient, on the other words, Ibragimov's coefficient, $ \  \beta(F,G) \  $ is defined by the formula

$$
\beta =\beta(F,G) \stackrel{def}{=}  \sup_{ A \in F, \  B \in G, \ {\bf P}(A) > 0 }  | {\bf P}(B/A) - {\bf P}(B)  |. \eqno(1.2)
$$
 \ Denote as usually  here and in the sequel by $ \   |\xi|_p \ $ the Lebesgue - Riesz $ \  L(p) \ $ norm of the r.v. $ \ \xi: \ $

$$
|\xi|_p = \left[  {\bf E} |\xi|^p  \right]^{1/p} :=
\left[  \int_{\Omega} |\xi(\omega)| \ {\bf P}(d \omega)  \right]^{1/p},  \ 1 \le p < \infty;
$$

$$
|\xi|_{\infty} := \vraisup_{\omega \in \Omega} |\xi(\omega)|.
$$

 \ Let $ \  \xi \in \in F, \ \xi \in L(p), \ \eta \in \in  G, \eta \in L(q), \ p,q \in [1, \infty]. \ $  Yu.A.Davydov in  [6]
proved  the following important inequality

$$
|\Cov(\xi,\eta)| \le 12 \ \alpha^{1 - 1/p - 1/q} \ |\xi|_p \ |\eta|_q, \ \frac{1}{p} + \frac{1}{q} < 1. \eqno(1.3)
$$
 \ The similar inequality for strong mixing coefficient $ \  \beta = \beta(F,G):  \  $

$$
|\Cov(\xi,\eta)| \le 2 \ \beta^{ 1/p} \ |\xi|_p \ |\eta|_q, \ \frac{1}{p} + \frac{1}{q} = 1, \eqno(1.4)
$$
or equally

$$
q = q(p) = \frac{p}{p-1} = p', \ p > 1; \ q(\infty) = 1,
$$
 may be found in the famous monograph of I.A.Ibragimov and Yu.A.Linnik [9]; see also the recent survey [19] and the article [12]. \par

  \ The following estimate  based only on the H\''older's inequality may be considered as trivial:

$$
|\Cov(\xi,\eta)| \le 2 \ |\xi|_p \ |\eta|_q, \ \frac{1}{p} + \frac{1}{q} = 1. \eqno(1.4a)
$$

 \ It is natural to expect that if the considered r.v. $ \  \xi, \ \eta \  $ have more light tails, for instance, satisfy the Kramer's condition,
then both the estimates (1.3) and (1.4) can be essentially strengthened. \par

 \ Both the inequalities (1.3) and (1.4) are essentially, i.e. up to multiplicative constants, non - refinable for all the  admissible
values $ \  p,q; \  $ see e.g. [6], for instance, on the classical probability space $ \  [0,1] \ $ equpped with Lebesgue measure there
exist two sigma - fields $ \ F, \ G $  and two non - zero symmetrical distributed r.v. $  \  \xi, \eta \  $ wor which

$$
|\Cov(\xi, \eta) | \ge \alpha(F,G)^{1 - 1/p - 1/q} \ |\xi|_p  \ |\eta|_q.
$$
 \ It is sufficient to take $ \ F = G. \  $\par

 \vspace{4mm}

 \   The   inequalities (1.3) and (1.4) are used in the investigation of the CLT for the dependent r.v., for the obtaining of
the non - asymptotical estimation for sums of  these r.v., in the statistics and  in the Monte - Carlo method, see e.g. [3], [4],  [12], [13] etc. \par

\vspace{4mm}

 {\bf  Our target in this short repor is extension of the estimates (1.3), (1.4) into the r.v. belonging to the
so - called Grand Lebesgue Spaces (GLS), in particular, into the exponential Orlisz spaces. }\par

\vspace{4mm}

  \  A modern result in this direction is represented in the article of E.Rio [17]; we intend to give  these covariation estimates
in another terms. \par

\vspace{4mm}

\section{\ Grand Lebesgue Spaces (GLS). Fundamental functions. }

 \vspace{4mm}

 \ Let $  Z = (\Omega, \cal{B},  {\bf P} ) $ be again the source probability space with non - trivial  normed  measure $ \ {\bf P}. \  $
 \ Let also  $  \psi = \psi(p), \ p \in [1, b), \ b = \const \in (1,\infty]  $ (or   $ p \in [1,b] $ ) be certain bounded
from below:  $  \ \inf \psi(p)  > 0 $ continuous inside the  {\it  semi - open} interval $   \ p \in [1, b) $ numerical function. We can and will suppose
$   \ b = \sup \{p, \psi(p) < \infty\},  \ $ so that  $ \  \supp \ \psi = [1, b) \  $  or $ \ \supp \ \psi = [1, b]. \ $ The set of all such a functions will be denoted by
$ \ \Psi(b) = \{  \psi(\cdot)  \}; \ \Psi := \Psi(\infty).  $\par

\ We agree to extend the definition these functions.  Indeed, we define for  arbitrary  $ \  \psi(\cdot) \in  \Psi(b) \  $  function in the case when
$  \ b < \infty \ $ for the values $ \ p > b \ $ formally as follows:

$$
\forall p > b \ \Rightarrow \psi(p) := \infty.
$$

 \ By definition, the (Banach) Grand Lebesgue Space  \ (GLS)  \ space   $  \ G\psi = G\psi(b)  $ consists on all the  (real or complex)  numerical  valued
measurable functions (random variables) $ \zeta  $ defined on our probability space $  \  \Omega \  $  and having a finite norm

$$
||\zeta|| = ||\zeta||G\psi \stackrel{def}{=} \sup_{p \in [1,b)} \left\{ \frac{|\zeta|_p}{\psi(p)} \right\}. \eqno(2.0)
$$

\vspace{4mm}

 \ These spaces  are Banach functional space, are complete, and rearrangement invariant in the classical sense, see [1],  chapters 1, 2.
 They were  investigated in particular in  many  works, see e.g.  [5], [7], [8], [10],  [11], [13], [14]. The function $ \  \psi = \psi(p) \ $
is said to be  {\it  the generating function } for this space. \par

 \ We refer here  some  used in the sequel facts about these spaces  and supplement more. \par

\vspace{4mm}

\  It is known that if  $  \ \zeta \ne 0, \ \zeta \in G\psi, \ $ then

  $$
  {\bf P}(|\zeta| > y) \le 2 \ \exp \left( \ - v_{\psi}^* (\ln (y/||\zeta||) \ \right), \ y \ge e \cdot ||\zeta||, \eqno(2.2)
  $$
 where $ \ v(p) = v_{\psi}(p) := p \ln \psi(p) \ $ and $  \  v^*(\cdot) \ $ denotes the Young - Fenchel, or Legendre  transform for
the function $ \ v(\cdot): $

$$
v^*(x) = \sup_{p \in \Dom (v)} (px - v(p)).
$$

  \ Conversely,  the last inequality may be reversed  in the following version: if

$$
{\bf P}(|\zeta| > y) \le 2 \ \exp \left(-v_{\psi}^* (\ln (y/K) \right). \ y \ge e \cdot K, \ K = \const \in (0,\infty),
$$
and if the function $ v_{\psi}(p), \ 1 \le p < \infty  \ $  is positive, finite for all the values $  \  p \in [1, \infty), $ continuous, convex
and such that

$$
\lim_{p \to \infty} \ln \psi(p) = \infty,
$$
then  $ \zeta \in G\psi $ and besides $  ||\zeta||  \le C(\psi) \cdot K:  $

$$
||\zeta||G\psi \le C_1 ||\zeta||L(M)  \le C_2 ||\zeta||G\psi, \ 0 < C_1 < C_2 < \infty. \eqno(2.3)
$$

 \ Moreover,  let us introduce the following  {\it exponential}  Young - Orlicz function

$$
N(u) = N_{\psi}(u) :=  \exp \left(v_{\psi}^* (\ln |u|) \right), \ |u| \ge e; \  N_{\psi}(u) = C u^2, \ |u| < e.
$$

 \ The Orlicz's norm $ \ ||\zeta||L(N_{\psi} ) \ $  is quite equivalent under formulated above conditions on the function  $ \ \psi(\cdot) \ $
to the GLS one:

$$
||\zeta||G\psi \le C_3 ||\zeta||L(N_{\psi})  \le C_4 ||\zeta||G\psi, \ 0 < C_1 < C_2 < \infty.
$$

 \  Furthermore,  let now $  \eta = \eta(z), \ z \in W $ be arbitrary family  of random variables  defined on any set $ W  $ such that

$$
\exists b = \const\in (1,\infty], \ \forall p \in [1,b)  \ \Rightarrow  \psi_W(p) := \sup_{z \in W} |\eta(z)|_p  < \infty. \eqno(2.4)
$$
 \ The function $  p \to \psi_W(p)  $ is named as a {\it  natural} function for the  family  of random variables $  W.  $  Obviously,

$$
\sup_{z \in W} ||\eta(z)||G\Psi_W = 1.
$$

 \ The family $ \ W \ $ may consists on the unique r.v., say $  \  \Delta: \ $

$$
\psi_{\Delta}(p):= |\Delta|_p,
$$
if of course  the last function is finite for some value $ \  p = p_0 > 1. \  $\par
  \ Note that the last condition is satisfied if for instance the r.v. $ \  \zeta \ $ satisfies the so - called Kramer's
condition; the inverse proposition is not true. \par

 \ Let us  bring two examples.  Define as usually the tail function for arbitrary numerical valued random variable $ \ \xi \ $

$$
T_{\xi}(y) \stackrel{def}{=} \max ( {\bf P}(\xi \ge y), \ {\bf P}(\xi \le -y) ), \ y \ge 0.
$$

{\bf Example 2.1.} Let $ \ m = \const > 0; $ define the function

$$
\psi_m(p) = p^{1/m}, \ p \in[1, \infty).
$$

 \ The tail inequality

$$
T_{\xi}(y) \le \exp \left(  - C y^m   \right), \ y \ge 0
$$
for some positive conctant $\ C \ $   is quite equvalent to the inclusion $ \   \xi \in G\psi_m.  \  $\par

\vspace{4mm}

{\bf Example 2.2.} Let $ \ b = \const > 1; \  \beta = \const > 0. $ Define the  following tail function

$$
T [b,\beta](y) : =   C \ y^{-b} \ (\ln y)^{\beta b - 1}, \ y \ge e,
$$
and the following $ \ \Psi(b) \ $ function with finite support

$$
\psi[b,\beta](p) =  (b-p)^{-\beta}, \ p \in [1,b); \  \psi[b,\beta](p)= \infty, \ p \ge b.
$$

 \ The  tail inequality of the form

$$
T_{\eta}(y) \le T[b,\beta](y), \ y \ge e
$$
entails  the inclusion $ \  \eta \in G\psi[b,\beta].  \ $  \par

 \ Note that the inverse proposition is not true. \par

\vspace{4mm}

 \ {\bf   Definition 2.1. }  The {\it fundamental function} for GLS $ \ G\psi_b \hspace{4mm} \phi[G\psi](\delta), \ \delta \in (0,\infty) $
may be calculated in accordance by the general theory of rearrangement invariant spaces [1], chapters 1,2  by a formula

$$
\phi[G\psi](\delta)  := \sup_{p \in [1,b)} \left\{  \frac{\delta^{1/p}}{ \psi(p) }  \right\}. \eqno(2.5)
$$
 \ This notion play a very important role in the Functional Analysis, theory of Fourier series, Operator Theory, Theory of Random Processes
 etc., see  the classical monograph [1]. For the GLS  this function was investigated in the preprint  [15]. It in proved in particular that
there exists a bilateral  continuous interrelation between fundamental  and generating function  for these spaces.\par

\vspace{4mm}

 \ {\bf   Definition 2.2. }  (See [15].) The {\it low truncated fundamental function} for the GLS $ \ G\psi_b, $ \hspace{3mm} namely,
$ \ \phi_s[G\psi](\delta), \ \delta \in (0,\infty), 0 < s < b \ $  is defined  by a formula

$$
\phi_s[G\psi](\delta)  := \sup_{p \in [s,b)} \left\{  \frac{\delta^{1/p}}{ \psi(p) }  \right\}, \ 1 \le s < b. \eqno(2.5a)
$$

\vspace{4mm}

 \ {\bf   Definition 2.3. }  (See [15].) The {\it upper truncated fundamental function} for the GLS
$ \ G\psi_b, $ \hspace{4mm}  indeed: $ \ \phi^s[G\psi](\delta), \ \delta \in (0,\infty), 0 < s < b \ $ is defined  by a formula

$$
\phi^s[G\psi](\delta)  := \sup_{p \in [s,b)} \left\{  \frac{\delta^{1/p}}{ \psi(p) }  \right\}, 1 \le s < b. \eqno(2.5b)
$$

\vspace{4mm}

 \ {\bf Example 2.1.a.} Let $ \ m = \const > 0; $ the fundamental function for the $ \  G\psi_m \  $ has a form

$$
\phi[G\psi_m](\delta) = (e m)^{-1/m} \ |\ln \delta|^{-1/m}, \ \delta \in (0, 1/e). \eqno(2.6)
$$

\vspace{4mm}

 \ {\bf Example 2.2.a}.  \ Define the following $ \  \Psi \ - \  $ function with finite support

$$
\tau_{b, \beta}(p) \stackrel{def}{=} (b-p)^{-\beta}, \ p \in [1,b). \eqno(2.7)
$$
 \ Here $ \ b = \const \in (1, \infty),\ \beta = \const \ge 0. \ $  The fundamental function for these space has a form

$$
\phi \left[G\tau_{b, \beta} \right](\delta) = \frac{b^{2 \beta -1} \ \beta^{\beta} \ \delta^{1/b}}{ |\ln \delta|^{\beta}}
=: K(b, \beta) \ \delta^{1/b} \ |\ln \delta|^{-\beta},  \ \delta \in (0, 1/e]. \eqno(2.8)
$$

\vspace{4mm}

\section{ Main results.  Strong mixing.}

 \vspace{4mm}

 \  Suppose  $  \  \xi \in G\psi, \ \eta \in G\nu  \ $ for certain $  \ \Psi \  $  functions $  \psi, \ \nu, $  and that

$$
  \xi \in \in F,  \ \eta \in \in  G. \eqno(3.0)
$$

 \ One can allow without loss of generality $ \   ||\xi||G\psi = ||\eta||G\nu = 1.  \ $ Then

$$
|\xi|_p \le \psi(p), \ | \eta|_{p/(p-1)} \le \nu(p/(p-1)), \ p \in [1,b). \eqno(3.1)
$$

 \ Define a new $ \ \Psi \ $ function

$$
\zeta(p) = \zeta[\psi,\nu](p) :=\psi(p) \ \nu(p/(p-1)); \eqno(3.2)
$$
then we have using the estimate (1.4)

$$
0.5 \ |\Cov(\xi, \eta)| \le \beta^{1/p}(F,G) \ \psi(p) \ \nu(p/(p-1)) = \beta^{1/p} \zeta[\psi,\nu](p),
$$
 therefore

$$
0.5 \ |\Cov(\xi, \eta)| \le \inf_p \left[ \beta^{1/p} \zeta(p) \right] =
$$

$$
\left\{ \sup_p \left[  \frac{\beta^{-1/p}}{\zeta[\psi,\nu](p)} \right]  \right\}^{-1} =
\frac{1}{\phi[G\zeta](1/\beta)}.
$$

 \ To summarize: \\

\vspace{4mm}

{\bf Theorem 3.1.} We  deduce under formulated above  notations and conditions

$$
|\Cov(\xi,\eta)|| \le  \frac{ 2 \ ||\xi||G\psi \ ||\eta||G\nu}{\phi[G\zeta[\psi,\nu]](1/\beta(F,G))}. \eqno(3.3)
$$

\vspace{4mm}

 \ Let us consider a particular case. \par

\vspace{4mm}

 \ {\bf  Definition 3.1.}  Let the function $ \ \psi = \psi(p) \ $ be from the set $ \  \Psi  = \Psi(\infty):  \ \supp \ \psi = [1, \infty). \ $
The function  $ \   \hat{\psi}  \  $ from this set is said to be {\it dual} to the function $  \ \psi(\cdot),  \ $ iff

$$
\hat{\psi}(p/(p-1)) = \psi(p); \ \Longleftrightarrow \hat{\psi}(p) = \psi(p/(p-1)). \eqno(3.4)
$$

\vspace{4mm}

 \ Evidently, $ \  \hat{\hat{\psi}} = \psi.  \  $ \par

\vspace{4mm}

 \ {\bf Proposition 3.1.} Suppose in addition to the conditions of theorem 3.1 that in (3.2)  $ \ \nu = \hat{\psi};  $ then
$ \  \zeta(p) = \psi^2(p) \ $ and following

$$
|\Cov(\xi,\eta)|| \le   2 \left[  \phi[G\psi] \left(\beta^{-1/2} \right) \right]^{-2}  \cdot ||\xi||G\psi \cdot ||\eta||G\hat{\psi}. \eqno(3.5)
$$

\vspace{4mm}

 \  {\bf Remark 3.1.} Let $ \ \nu(\cdot) \in \Psi(b), \ b = \const \in (1, \infty].  $ We assert that the Grand Lebesgue Space
$ \  G\hat{\nu} \ $  consists only on the essentially bounded variables:

$$
||\zeta||  G\hat{\nu} \asymp |\zeta|_{\infty}. \eqno(3.6)
$$

\vspace{4mm}

 \ {\bf Proof.} The inclusion $ \ L_{\infty} \subset  G\hat{\nu} \ $  is evident; we must ground  an inverse inclusion. \par
 \ So,  let  $  \ \zeta \in  G\hat{\nu}, \ \supp (\nu) = [1, b). \ $ We have taking into account the continuity of the function $ \  \nu(\cdot) \ $
at the point $ \ 1 + 0:  $

$$
\lim_{p \to \infty} \hat{\nu}(p) = \lim_{p \to \infty} \nu \left( \frac{p}{p-1}    \right) = \nu(1) < \infty,
$$
therefore

$$
 \overline{\lim}_{p \to \infty} |\zeta|_p < \infty \ \Longleftrightarrow \vraisup_{\omega \in \Omega}|\zeta(\omega)| < \infty.
$$

\vspace{4mm}

\section{ Main results.  Uniform  (Rosenbatt) mixing.}

 \vspace{4mm}

 \ Let as before $ \ \ \psi, \ \nu \ $ be two $ \ \Psi \ $ functions and  let $  \ \alpha, \ \beta = \const \in [0,1]. \ $

 \ Denote by $ \ T \ $ the domain in the  positive quarter plane

$$
T = \{ p,q: \  p,q \ge 1, 1/p + 1/q < 1 \};
$$
 ``T'' implies a triangle for the inverse values $  \  x = 1/p, \ y = 1/q.  \  $  Introduse the following functions

$$
\Phi[\psi, \nu](\alpha, \beta) \stackrel{def}{=}
\sup_{(p,q) \in T} \left[ \frac{\alpha^{1/p} \ \beta^{1/q}}{\psi(p) \ \nu(q)}   \right],  \eqno(4.1)
$$

$$
\theta[\nu]_{\beta}(p) := \frac{\psi(p)}{\phi^{(p')}[G\nu](\beta)}, \eqno(4.2)
$$

 \ We have

$$
\Phi[\psi, \nu](\alpha, \beta) = \sup_p \left\{   \frac{\alpha^{1/p}}{\psi(p)} \ \sup_{q \ge p'} \frac{\beta^{1/q}}{ \nu(q)}   \right\} =
$$

$$
\sup_p \left\{   \frac{\alpha^{1/p}}{\psi(p)} \ \sup_{q \ge p'} \frac{\beta^{1/q}}{ \nu(q)}   \right\} =
\sup_p \left[   \frac{\alpha^{1/p}}{ \theta[\nu]_{\beta}(p)}   \right] =
\phi[G\theta_{\beta} ](\alpha). \eqno(4.3)
$$

 \vspace{4mm}

 \ {\bf  Theorem 4.1.}

$$
|\Cov(\xi, \ \eta)| \le  \frac{12 \ \alpha \ ||\xi||G\psi \ ||\eta||G\nu } {\Phi[\psi, \nu](\alpha, \beta)} =
 \frac{12 \ \alpha \ ||\xi||G\psi \ ||\eta||G\nu }{\phi[G\theta_{\alpha} ](\alpha)}.\eqno(4.4)
$$
 \ Here $ \   \alpha = \alpha(F,G). \ $ \par

\vspace{4mm}

 \ {\bf Proof.} Assume $  \  ||\xi||G\psi = ||\eta||G\nu = 1. \  $ Then

$$
|\xi|_p  \le \psi(p), \ |\eta|_q \le \nu(q), \   (p,q) \in T.
$$

 \ One can apply the Davydov's inequality (1.3):

$$
(12 \alpha)^{-1} |\Cov(\xi, \ \eta)| \le \alpha^{-1/p} \ \alpha^{-1/q} \ \psi(p) \ \nu(q),
$$
therefore

$$
(12 \alpha)^{-1} |\Cov(\xi, \ \eta)| \le  \inf_{  (p,q) \in D }  \left[\alpha^{-1/p} \ \alpha^{-1/q} \ \psi(p) \ \nu(q) \right] =
$$

$$
\frac{ 1 }{\phi[G\theta_{\alpha} ](\alpha)}  =  \frac{ ||\xi||G\psi \ ||\eta||G\nu }{\phi[G\theta_{\alpha} ](\alpha)}, \eqno(4.5)
$$
Q.E.D.

 \vspace{4mm}

 \  Let's turn again our attention  on the function $ \ \Phi[\psi, \nu](\alpha, \beta) \ $  from the definition (4.1). In all the considered
examples it  allows a {\it factorization }

$$
\Phi[\psi, \nu](\alpha, \beta) = \phi[G\psi](\alpha) \cdot \phi[G\nu](\beta) \eqno(4.6)
$$
for all sufficiently small values $  \  \alpha, \ \beta. \  $ We intend further to investigate the possibility of this relation (4.6). \par

 \vspace{4mm}

  \ We need first of all to return  to the investigation of the fundamental function for GLS. \par

\vspace{4mm}

  \ {\bf Case A. Infinite support.} \par

\ Let  $ \psi(\cdot)  \in \Psi(\infty) = \Psi.\ $  Denote  for an arbitrary $ \ \Psi \ $ function the following transform

$$
g(x) = g[\psi](x):=  -\ln \psi(1/x), \ x \in (0,1), \eqno(4.7)
$$

$$
p_0 = p_0(\delta) := \argmax_{p \in [1, \infty)} \left(  \frac{\delta^{1/p}}{\psi(p)}  \right), \eqno(4.8)
$$
so that $ \ x_0 = x_0(\delta) = 1/p_0(\delta)  \ $ and

$$
  \frac{\delta^{1/p_0}}{\psi(p_0)}  = \phi[G\psi](\delta).
$$

\vspace{4mm}

\ {\bf Lemma A.}  Suppose that the derivative $  \ g'(x) = g[\psi]'(x) $ there exists, is continuous  in the open interval (0,1),  is
strictly decreasing and such that

$$
\lim_{x \to 0} g'[\psi](x) = \infty. \eqno(4.9)
$$
 Then

$$
\lim_{\delta \to 0+}  p_0(\delta)  = \infty.  \eqno(4.10)
$$

\vspace{4mm}

 \ {\bf Proof} follows immediately from the equation

$$
 g' [\psi](1/p_0(\delta) )=   g'(x_0) = g[\psi]'(x_0)  = \ln (1/\delta), \ \delta \in (0,1). \eqno(4.11)
$$

\vspace{4mm}

  \ {\bf Case B. Finite support.} \par

\ Let now  $ \psi(\cdot)  \in \Psi(b), \ $ where $ \  1 < b = \const < \infty.  \  $  Denote as before

$$
g(x) = g[\psi](x):=  -\ln \psi(1/x), \ x \in (0,1/b),
$$

$$
q_0 = q_0(\delta) := \argmax_{q \in [1, b)} \left(  \frac{\delta^{1/q}}{\psi(q)}  \right). \eqno(4.12)
$$

\vspace{4mm}

\ {\bf Lemma B.}  Suppose that the derivative $  \ g'(x) = g[\psi]'(x) $ there exists, is continuous  in the open interval $ \ (0,1/b), \ $  is
strictly increasing and such that

$$
\lim_{x \to 1/b} g'(x) = \infty. \eqno(4.13)
$$
 Then  as $ \ \delta \downarrow 0+    \ $

$$
  q_0(\delta) \uparrow  1/b.  \eqno(4.14)
$$

\vspace{4mm}

 \ {\bf Proof} follows immediately likewise before from the equation  (4.11). \par

\vspace{4mm}

 \ Let us return to the factorization equality (4.6).  Introduce the following  ``rectangle''

$$
R := \{  (p,q): 1 \le p,q < \infty \}.
$$

 \ Obviously,

$$
\Phi[\psi, \nu](\alpha, \beta) =
\sup_{(p,q) \in T} \left[ \frac{\alpha^{1/p} \ \beta^{1/q}}{\psi(p) \ \nu(q)}   \right] \le
$$

$$
\sup_{(p,q) \in R} \left[ \frac{\alpha^{1/p} \ \beta^{1/q}}{\psi(p) \ \nu(q)}   \right]  =
\phi[G\psi](\alpha) \cdot \phi[G\nu](\beta). \eqno(4.14)
$$
 \ It  remains to ground the opposite inequality, of course, for sufficiently smallest values $ \ \alpha \ $ and $ \ \beta. \ $  \par

\vspace{4mm}

 \ We consider consequently three cases. \par

\vspace{4mm}

\ {\bf Case 4.1. Infinite supports.} \par

 \ Suppose the two functions  $  \ \psi = \psi(p), \ \nu = \nu(p) $ belonging to the set $ \ \Psi(\infty) = \Psi \ $ be a given and
both these functions satisfy to the conditions of lemma A,  in particular, the relation  (4.9):

$$
\lim_{x \to 0} g'[\psi](x) = \lim_{x \to 0} g'[\nu](x) = \infty.
$$

 \ We  have using Lemma A:

$$
\phi[G\psi](\alpha) \cdot \phi[G\nu](\beta)  =
 \frac{\alpha^{1/p_0(\alpha)}}{\psi_m(p_0(\alpha))}   \cdot \frac{\beta^{1/p_0(\beta)}}{\psi_n(p_0(\beta))} , \eqno(4.15)
$$
as long as

$$
\lim_{\alpha \to 0+} p_0(\alpha) =  \lim_{\beta \to 0+} p_0(\beta) = \infty.
$$

 We deduce therefore for sufficiently smallest values $ \ \alpha \ $ and $ \ \beta \ $

$$
\frac{1}{p_0(\alpha)} +  \frac{1}{p_0(\beta)} < 1,
$$
 on the other words   the optimal pair  $  \ ( p_0(\alpha), p_0(\beta)  )  \  $ belongs to the set $ \ T. \ $ \par

  \ To be more specifically, note that if for instance

$$
\alpha_0  := \exp \left( - g'[\psi](1/e)  \right)
$$
and correspondingly

$$
\beta_0  := \exp \left( - g'[\nu](1/e)  \right),
$$
 we deduce under asumptions of lemma A  that when $ \ \alpha \le \alpha_0, \ \beta \le \beta_0   \ $

$$
p_0(\alpha) \ge e, \ p_0(\beta) \ge e
$$
whence

$$
\frac{1}{p_0(\alpha)}  + \frac{1}{p_0(\beta)} < 1.
$$

\vspace{4mm}

 \ {\bf Case 4.2. Finite supports.} \par

  \ Let two functions $ \ \psi = \psi(p) \ $ and $ \   \nu =\nu(p) \ $ be a given;  suppose that these
functions belongs correspondingly  to the sets

$$
\psi(\cdot) \in G\psi(b_1), \  \nu(\cdot) \in G\psi(b_2), \ b_1, b_2 > 1.
$$

  \ Assume  that the derivatives $  \  g[\psi]'(x) $  and  $  \ g[\nu]'(x) $ there exist, are continuous  in the open
intervals $ \ (0,1/b_{1,2}), \ $  are strictly increasing and such that

$$
\lim_{x \to 1/b_1} g'[\psi](x) = \infty =  \lim_{x \to 1/b_2} g'[\nu](x).
$$

 \ We conclude using Lemma B that the relation (4.6) there holds if

$$
\frac{1}{b_1} +  \frac{1}{b_2} < 1. \eqno(4.16)
$$

\vspace{4mm}

 \ To be more concrete, suppose estimate (4.16) be satisfied.  Define the following values

$$
\beta_1 = \exp \left(  - g'[\psi]  \left( \frac{b_1 + 1}{3b_1}  \right)  \right),
$$

$$
\beta_2 = \exp \left(  - g'[\nu]  \left( \frac{b_2 + 1}{3b_2}  \right)  \right),
$$
or equally

$$
q(\beta_j) = \frac{3b_j}{b_j+1}, \ j = 1,2.
$$
  \ We conclude for  all the values $ \  \theta_j \in (0, \beta_j) \ $

$$
\frac{1}{q(\theta_1)} + \frac{1}{q(\theta_2)} \le  \frac{1}{q(\beta_1)} + \frac{1}{q(\beta_2)} =
$$

$$
\frac{b_1 + 1}{3 b_1} + \frac{b_2 + 1}{3 b_2} = \frac{2}{3} + \frac{1}{3b_1} + \frac{1}{3b_2} < 1.
$$

\vspace{4mm}

 \ {\bf Case 4.3. ``Mixed `` \ case.} \par

 \ Assume   $   \  \psi(\cdot) \in  \ \Psi(\infty),  \ \nu(\cdot) \in G\psi_{b}, \ b = \const > 1,  $ and

$$
\lim_{x \to 0} g'[\psi](x) = \infty =  \lim_{x \to 1/b} g'[\nu](x).
$$
 \ We conclude likewise foregoing propositions that  then the equality  (4.6) there holds still in this case, since

$$
 \frac{1}{\infty} + \frac{1}{b} = \frac{1}{b} < 1.
$$

 \ To be more precisely,  it is sufficient to pick the (positive) values  $ \ \alpha_0 \ $ and $  \ \beta_0 \ $ as follows:
$ \ \alpha_0, \ \beta_0  \in (0,1) \ $ and

$$
g'[\psi] \left(  \frac{b-1}{3b}  \right) = |\ln \alpha_0|,
$$

$$
g'[\nu]  \left(  \frac{b+1}{2b}   \right) =  |\ln \beta_0|.
$$

 \ Then $ \  p_0 := \frac{3b}{b-1}  \ $ and $ \  q_0 := \frac{2b}{b+1}, \ $  so that $ \ p_0 > 1, $

$$
 1 < q_0  = \frac{2b}{b+1} < b,
$$
and

$$
\frac{1}{p_0}   +  \frac{1}{q_0} = \frac{b-1}{3b} + \frac{b+1}{2b} <  \frac{b-1}{2b} + \frac{b+1}{2b} = 1,
$$
as long as $  \ b > 1. $\par

\vspace{4mm}

\section{ Main results. The case of identical spaces. Examples.}

 \vspace{4mm}

 \ We suppose in addition to the conditions (and notations) of the last section that the $ \ \Psi \ $  functions
$ \ \psi(\cdot), \ \nu(\cdot)  $ coinsides: $ \ \psi(\cdot) = \nu(\cdot).  $\par

 \ In detail: assume $  \xi \in \in F, \ \xi \in G\psi, \ \eta \in \in G, \ \eta \in G\psi.  $ \par

\vspace{4mm}

{\bf Theorem 5.1.}

$$
|\Cov(\xi,  \ \eta)|  \le 12 \ \alpha(F,G) \cdot
\frac{ ||\xi||G\psi \ ||\eta||G\psi}{\phi^2[G\psi](\alpha(F,G))}. \eqno(5.1)
$$

\vspace{4mm}

 \ {\bf Proof.}  \ Suppose for simplicity $  \ ||\xi||G\psi = 1 = ||\eta||G\psi; $ then

$$
|\xi|_p  \le \psi(p), \ |\eta|_p \le \psi(p), \   p \in [1,b),  \  b = \const \in (1, \infty].
$$

 \ We can apply again the Davydov's inequality (1.3):

$$
(12 \alpha)^{-1} |\Cov(\xi, \ \eta)| \le \alpha^{-2/p} \  \psi^2(p),
$$
therefore

$$
(12 \alpha)^{-1} |\Cov(\xi, \ \eta)| \le \inf_p \left[ \alpha^{-2/p} \  \psi^2(p) \right]  =
$$

$$
 \left[ \  \sup_p \frac{\alpha^{1/p}}{\psi(p)} \right]^{-2} =  \left[ \  \phi[G\psi](\alpha) \right]^{-2} =
\frac{ ||\xi||G\psi \ ||\eta||G\psi}{\phi^2[G\psi](\alpha)}, \eqno(5.2)
$$
Q.E.D.

\vspace{4mm}

 \ {\bf Example 5.1. Infinite supports.} \par

 \  Define the following $ \ \Psi \ $ functions

$$
\psi_m(p) = p^{1/m}, \ \psi_n(p) = p^{1/n}, \ p \in [1, \infty), \ m,n = \const > 0.
$$
 \ If $ \xi \in \in F, \ \xi \in G\psi_m, \ \eta \in \in G, \ \eta \in G\psi_n \ $ and $  \   \alpha = \alpha(G,F) \le e^{-1},  \  $ then

$$
|\Cov(\xi, \eta) \le 12 \ e^{1/m + 1/n} \ m^{1/m} \ n^{1/n} \ \alpha \ |\ln \alpha|^{1/m + 1/n} \ ||\xi||G\psi_m \ ||\eta||G\psi_n.
\eqno(5.3)
$$

\vspace{4mm}

 \ {\bf Example 5.2. Finite supports.} \par

 \ Recall the definition of  the following $ \  \Psi \ - \  $ function

$$
\tau_{b, \beta}(p) \stackrel{def}{=} (b-p)^{-\beta}, \ p \in [1,b). \eqno(5.4)
$$
 \ Here $ \ b = \const \in (1, \infty),\ \beta = \const \ge 0. \ $  As we knew, the fundamental function for these space has a form

$$
\phi \left[G\tau_{b, \beta} \right](\delta) = \frac{b^{2 \beta -1} \ \beta^{\beta} \ \delta^{1/b}}{ |\ln \delta|^{\beta}}
=: K(b, \beta) \ \delta^{1/b} \ |\ln \delta|^{-\beta},  \ \delta \to 0+. \eqno(5.5)
$$

 \ This relation allow us to calculate the required covariation. Namely,
 if $ \xi \in \in F, \ \xi \in G\tau_{b_1, \beta_1}, \ \eta \in \in G, \ \eta \in G\tau_{b_2, \beta_2}, \
b_{1,2} = \const \in (1, \infty),\ \beta_{1,2} = \const \ge 0 \ $ and $  \   \alpha = \alpha(G,F) \le e^{-1},  \  $ then

$$
|\Cov(\xi, \eta)| \le 12 \ K(b_1, \beta_1) \ K(b_2, \beta_2) \ \alpha^{1 - 1/b_1 - 1/b_2} \ |\ln \alpha|^{\beta_1 +\beta_2} \ \times
$$

$$
||\xi||G\tau_{b_1, \beta_1} \ ||\eta||G\tau_{b_2, \beta_2},  \eqno(5.6)
$$
if of course $   \  1/b_1 + 1/b_2 < 1. \  $\par

\vspace{4mm}

 \ {\bf Example 5.3. ``Mixed `` \ case.} \par

\vspace{4mm}

 \ Assume  $ \ \xi \in \in F, \ \xi \in G\psi_m, \  m  = \const > 0; \ \eta \in \in G, \ \eta \in G\tau_{b, \beta},
b = \const > 1, \ \beta = \const >0, \ $ and denote as before
$ \  \alpha = \alpha(F,G).  \ $ We obtain after some calculations

$$
|\Cov(\xi, \eta)| \le 12 \  (em)^{1/m} \ K(b, \beta) \  \alpha^{1 - 1/b} \ |\ln \alpha|^{\beta +1/m} \ \times
$$

$$
||\xi||G\psi_m \ ||\eta||G\tau_{b, \beta}. \eqno(5.7)
$$

 \ {\bf Remark 5.1.}  For the ``greatest'' values  $ \ \alpha \ $ and $ \ \beta, \ $ say $ \ \alpha  \ge 1/e, $ one can use
the trivial estimate 1.4a. \par

\vspace{4mm}

 \ {\bf Example 5.4.  ``Combined'' event.} \par

\vspace{4mm}

 \ Suppose here that $   \   \xi \in G\psi_b \ $ for some $ \  b = \const   > 1 \ $ and that $ \  \eta \in L(q(0)),  \  $
where $  \  q(0)' < b. \  $  We derive consequently $ \ |\xi|_p \le ||\xi||G\psi \cdot\psi(p), \ 1 \le p < b;  $

$$
(12 \alpha)^{-1} \ \alpha^{1/q(0)} (|\Cov(\xi,\eta))| ||\xi||G\psi \ |\eta|_{q(0)} ) \le \alpha^{-1/p} \psi(p), \ p \ge q'(0),
$$
whence

$$
(12 \alpha)^{-1} \ \alpha^{1/q(0)} \ ( |\Cov(\xi,\eta)|) /( ||\xi||G\psi \ |\eta|_{q(0)} ) \le  =
$$

$$
  \inf_{ p \in [q'(0), b)  }\alpha^{-1/p} \psi(p)  = \frac{1}{\phi_{q_0}[G[\psi](\alpha)}.
$$
 \ Thus, we deduce in the considered case

$$
|\Cov(\xi,\eta)| \le 12 \ \alpha^{1 - 1/q(0)} \cdot \frac{ ||\xi||G\psi \cdot |\eta|_{q(0) }}{\phi_{q(0)} [G\psi](\alpha) }.
$$

%&

\vspace{4mm}

\section{ Application to the classical CLT.}

 \vspace{4mm}

 \ Let  $ \  \gamma(i), \ i = 0, \pm 1, \pm 2 \ldots  \ $ be a centered strictly  stationary sequence of  r.v.  A new denotations:

$$
S(n) := n^{-1/2} \sum_{i=1}^n \gamma(i), \ n = 1,2,\ldots;. \eqno(6.0)
$$

$$
\Sigma(n) : = \Var (S(n)),  \ \Sigma:= \lim_{n \to \infty} \Sigma(n);   \eqno(6.1)
$$

$$
F_0 := \sigma \{ \gamma(i), \ i \le 0 \};  \ F^k := \sigma \{ \gamma(j),  \ j \ge k \};
$$

$$
\alpha(k) := \alpha \left( F_0, F^k  \right); \ \beta(k) :=  \beta \left( F_0, F^k  \right). \eqno(6.2)
$$

 \ Introduce also the following $ \ \Psi \ $ function as a natural function for the sequence $ \ \{\gamma(i)\}: \ $

$$
\psi[\gamma](p) := |\gamma(0)|_p; \eqno(6.3)
$$
if of course there exists for {\it some} value $ \ p \ $ greatest than one. \par

\vspace{4mm}

 \ Recall that the sequence $ \ \gamma(\cdot) \ $ satisfies the CLT, iff $ \  \exists \ \Sigma \in (0,\infty)  \ $ and
the sequence $ \ \{ S(n) \} \ $ converges in distribution as $ \ n \to \infty \ $ to the centered normal (Gaussian) law
with variance $ \ \Sigma. \ $ \par
 \ The methods of obtaining CLT for the random  sequences  satisfying some mixing conditions are well known,  see e.g.  [9], [16], [18],
 [19] etc.  The essential  moment for this proof is the following variation estimate, which follows immediately from the
foregoing covariation estimates. \par

\vspace{4mm}

{\bf Theorem 6.1.} Assume that the function  $ \ \psi[\gamma](p) \ $ is not trivial: $ \ \exists p_0 > 1 \ \Rightarrow
 \psi[\gamma](p_0) < \infty. \ $ Define also  a numerical  positive sequence

$$
y(k) = y[\gamma](k) := \frac{ \alpha(k)}{\phi^2[G\psi[\gamma]](\alpha(k))}, \ k = 2,2,\ldots. \eqno(6.4)
$$
 \ If

$$
\sum_{k=2}^{\infty} y[\gamma](k) < \infty,
$$
then the value $  \Sigma \ $  there exists and is finite: $ \ \Sigma  \in [0,   \infty). $ \par

\vspace{4mm}

 \ Define now a new $ \ \Psi \ $ function

$$
\zeta[\psi](p) :=\psi(p) \ \psi(p/(p-1)). \eqno(6.5)
$$
and the correspondent numerical sequence, also positive

$$
 z(k) = z[\gamma](k) :=  \frac{1}{\phi[G\zeta[\psi]](1/\beta(k))}. \eqno(6.6)
$$

\vspace{4mm}

{\bf Theorem 6.2.} We  deduce under formulated above  notations and conditions that if

$$
\sum_{k=2}^{\infty} z[\gamma](k) < \infty,
$$

then the value $  \Sigma \ $  there exists and is finite: $ \ \Sigma  \in [0,   \infty). $ \par

\vspace{4mm}

\section{ \ Concluding remarks.}

 \vspace{4mm}

\ {\bf A.}  Offered here results may be easily generalized onto another types of mixing, as well as onto others
r.i. spaces: Lorentz,  Marcinkiewicz etc. instead  GLS.  All we need - the source  $ L(p), \ L(q) $ estimate of the form

$$
|\Cov(\xi, \ \eta)| \le h(p,q) \ |\xi|_p \ |\eta|_q, \ (p,q) \in D, \eqno(7.1)
$$
where $ \ D \ $ is arbitrary domain in the positive quarter plane. \par

 \ Suppose $ \ \xi \in G\psi, \  \eta \in G\nu. $  Then

$$
|\Cov(\xi, \ \eta)| \le \inf_{(p,q) \in D} \left[  h(p,q) \ \psi(p)\ \nu(q) \right] \ ||\xi||G\psi \ ||\eta||G\nu. \eqno(7.2)
$$

\vspace{4mm}

 \ {\bf B.} It is interest by our opinion to derive the lower bound for considered covariance for diffent Banach spaces. \par

\vspace{4mm}

\ {\bf C.}  Define the following $ \ \psi_{(r) } = \psi_{(r) } (p), \ r = \const > 1, \ - $ function as follows:
$ \  b(\psi_{(r) }) = r  \ $ and

$$
\psi_{(r) }(p) = 1, \ 1 \le p \le r.
$$
 One can define formally $ \psi_{(r) }(p) = \infty, \ p > r.  $ Then

$$
||f||G\psi_{(r) } =  |f|L_r.
$$
 Thus, the classical Lebesgue - Riesz spaces $  \ L_r \ $ are a particular, more precisely, extremal cases of the Grand Lebesgue ones. \par

 \ As long as the fundamental function for $ \   L_r(\Omega) \ $ space  builded by atomless measure $ \  {\bf P} \  $ is equal

$$
\phi[L_r]  (\delta) = \delta^{1/r}, \ \delta \in [0,1],
$$
 we derive the Davydov's estimate  (1.3)  in turn from the proposition of theorem 4.1; as well as  the inequality (1.4) follows from theorem
(3.1).\par

 \ As a slight consequence: both the assertions of theorems (3.1) and (4.1)  are  in general case  essentially non - improvable. \par

\vspace{4mm}

{\ \bf D.}  The case of the so - called   $ \ \rho \ - $ mixing is very simple for covariation estimation, see e.g. (20). \par

 \vspace{6mm}

 {\bf References.}

 \vspace{4mm}

{\bf 1. Bennet C., Sharpley R.}  {\it  Interpolation of operators.} Orlando, Academic
Press Inc., (1988). \\

\vspace{3mm}

{\bf 2. Billingsley P. }  {\it  Convergence of probability measures.} New York : Wiley 1968.\\

 \vspace{3mm}

{\bf 3. Bradley, R.C.}  (2007). {\it Introduction to strong mixing conditions. Vol. 1,2,3.}
Kendrick Press.

\vspace{3mm}

{\bf 4. Bryc, W. and Dembo, A. }  {\it Large deviations and strong mixing.}  Ann. Inst.
Henri Poincar'e, 32. (1996), pp. 549–569.

\vspace{3mm}

{\bf 5.  Buldygin V.V., Kozachenko Yu.V. }  {\it Metric Characterization of Random
Variables and Random Processes.} 1998, Translations of Mathematics Monograph, AMS, v.188. \\

\vspace{3mm}

{\bf 6. Davydov Yu.A.} {\it Convergence of distributions generated by stationary stochastic processes. }
 Theory Probab. Appli. 13, 691-696.\\

\vspace{3mm}

 {\bf 7. A. Fiorenza.}   {\it Duality and reflexivity in grand Lebesgue spaces. } Collect. Math.
{\bf 51,}  (2000), 131  - 148. \\

 \vspace{3mm}

{\bf  8. A. Fiorenza and G.E. Karadzhov.} {\it Grand and small Lebesgue spaces and
their analogs.} Consiglio Nationale Delle Ricerche, Instituto per le Applicazioni
del Calcoto Mauro Picone”, Sezione di Napoli, Rapporto tecnico 272/03, (2005).\\

\vspace{3mm}

{\bf 9. Ibragimov I.A., Linnik Yu.A.} {\it Independent and stationary sequences of
random variables.} Groningen Wolters-Noordhoff, 1971 \\

\vspace{3mm}

{\bf 10.  T. Iwaniec and C. Sbordone.} {\it On the integrability of the Jacobian under minimal
hypotheses. } Arch. Rat.Mech. Anal., 119, (1992), 129-143. \\

 \vspace{3mm}

{\bf 11. Kozachenko Yu. V., Ostrovsky E.I. }  (1985). {\it The Banach Spaces of random Variables of subgaussian Type. } Theory of Probab.
and Math. Stat. (in Russian). Kiev, KSU, 32, 43-57. \\

\vspace{3mm}

{\bf  12. Florence Merlev'ede, Magda Peligrad, Emmanuel Rio.} {\it Bernstein inequality and moderate deviations under
strong mixing conditions}.  Edited by Christian Houdr'e, Vladimir Koltchinskii,
David M. Mason and Magda Peligrad. High dimensional probability V : the 5th International
Conference (HDP V), May 2008, Luminy, France. Institute of Mathematical Statistics, Beachwood, OH, pp.273-292, 2009.
<inria-00360856>

 \vspace{3mm}

{\bf 13. Ostrovsky E.I. } (1999). {\it Exponential estimations for Random Fields and its
applications,} (in Russian). Moscow-Obninsk, OINPE. \\

 \vspace{3mm}

{\bf 14. Ostrovsky E. and Sirota L. } {\it   Entropy and Grand Lebesgue Spaces approach for the problem  of Prokhorov - Skorokhod
continuity of discontinuous random fields. }\\
arXiv:1512.01909v1 [math.Pr] 7 Dec 2015 \\

 \vspace{3mm}

{\bf 15. Ostrovsky E. and Sirota L. } {\it  Fundamental function for Grand Lebesgue Spaces.  }
arXiv:1509.03644v1  [math.FA]  11 Sep 2015 \\

\vspace{3mm}

{\bf  16.  Magda Peligrad and Sergey Utev.} {\it Central Limit Theorem for linear processes.}
The Annals of Probability, 1997, Vol 25, N 1, 443 - 456. \\

\vspace{3mm}

{\bf 17. Rio E.} {\it Covariance inequalities for strongly mixing processes.}   1993,
Ann. Inst. of H. Poincare, {\bf 29},  589 - 597.  \\

\vspace{3mm}

{\bf 18.  S. Utev, M. Peligrad.} {\it Maximal inequalities and an invariance principle for a class of weakly dependent random variables.}
Journal of Theoretical Probability, 16 (2003). 101 \ - \ 115.\\

\vspace{3mm}

{\bf 19. K.Yoshishara.}  {\it  Moment inequalities for mixing sequences.   } Kodai Math. J.,
1, (1978), 316-328.

\vspace{3mm}

{\bf 20. Q.Y. Wu, Y.Y. Jiang.} {\it Some strong limit theorems for   mixing sequences of random variables.} Statistics and   Probability Letters,
{\bf 78,} \ (2008), 1017 \ - \ 1023.

\vspace{3mm}

\end{document}